\documentstyle[11pt,amsfonts,fullpage]{amsart}

\newtheorem{theorem}{Theorem}[section]

\newtheorem{corollary}[theorem]{Corollary}

\newtheorem{lemma}[theorem]{Lemma}

\newtheorem{proposition}[theorem]{Proposition}

\def\endproof{\qed \medskip}
\def\blacksquare{\hbox to .60em{\vrule width .60em height .60em}}

\begin{document}

\title[]{On Uniqueness and Differentiability in the Space of Yamabe 
Metrics}

\author[]{Michael T. Anderson}

\thanks{Partially supported by NSF Grant DMS 0072591 and 0305865}


\maketitle

\abstract
It is shown that there is a unique Yamabe representative for a 
generic set of conformal classes in the space of metrics on any manifold. 
At such classes, the scalar curvature functional is shown to be 
differentiable on the space of Yamabe metrics. In addition, some 
sufficient conditions are given which imply that a Yamabe metric of 
locally maximal scalar curvature is necessarily Einstein. 
\endabstract

\setcounter{section}{0}

\section{Introduction.}
\setcounter{equation}{0}

 Let $M$ be a closed $n$-dimensional manifold. For a given smooth 
metric $g$ on $M$, let $[g]$ denote the conformal class of $g$, 
consisting of smooth metrics on $M$ pointwise conformal to $g$. By the 
solution to the Yamabe problem [1], [6], in each conformal class $[g]$ 
there is a Yamabe metric $\gamma$ of constant scalar curvature 
$s_{\gamma}$; the metric $\gamma$ minimizes the total scalar curvature 
(or Einstein-Hilbert action)
$${\cal S} (\bar g) = v^{-(n-2)/n}\int_{M}s_{\bar g}dV_{\bar g}. $$
when ${\cal S}$ is restricted to the class of conformal metrics $\bar g 
\in [g]$. Here $s_{g}$ denotes the scalar curvature of $g$, $dV_{g}$ 
the volume form, and $v$ the total volume of $(M, g)$. The sign of 
$s_{\gamma}$, i.e $s_{\gamma} < 0$, $s_{\gamma} = 0$ or $s_{\gamma} > 
0$ depends only on the conformal class $[\gamma]$.

 A number of general features of the class of Yamabe metrics of 
non-positive scalar curvature are well understood, cf. [2], [7] for 
example. Thus, negative conformal classes have a unique unit volume 
Yamabe metric. The space ${\cal Y}^{-}$ of all unit volume negative 
Yamabe metrics forms a smooth infinite dimensional manifold ${\cal 
Y}^{-}$, transverse to the space of conformal classes, in the space of 
all unit volume metrics on $M$. The scalar curvature $s$ defines a 
smooth function $s: {\cal Y}^{-} \rightarrow  {\Bbb R}$, whose critical 
points are exactly Einstein metrics (of negative scalar curvature) on 
$M$. Similar results hold for the space of non-positive Yamabe metrics. 
All of these results essentially derive from the fact that Yamabe 
metrics satisfy an elliptic equation, the Yamabe equation, whose 
solutions satisfy a maximum principle when $s_{\gamma} \leq 0$.

 It has been an open issue for some time to what extent such general 
features continue to hold for the space ${\cal Y}^{+}$ of positive unit 
volume Yamabe metrics on $M$ where the corresponding Yamabe equation 
does not satisfy a maximum principle. Thus, in general, uniqueness of 
Yamabe metrics fails for positive conformal classes. For example, the 
conformal class $[\gamma_{0}]$ of the round metric $\gamma_{0}$ on 
$S^{n}$ admits a large, non-compact, family of Yamabe metrics. On 
$S^{n-1}\times S^{1}$, there are 1-parameter families of Yamabe metrics 
in the conformal class of a product metric, for infinitely many values 
of radii of $S^{n-1}$ and $S^{1}$, cf. [7], [8].

 The purpose of this note is to establish some partial answers to these 
issues. Some of the main results are as follows:
\begin{theorem} \label{t 1.1.}
  Generically, Yamabe metrics are unique in their conformal class. 
Thus, there is an open and dense set ${\cal U}$ in the space of 
positive conformal classes such that each $[g]\in{\cal U}$ has a unique 
unit volume Yamabe metric $\gamma\in [g]$.
\end{theorem}

 We also prove that the scalar curvature function $s: {\cal Y}^{+} 
\rightarrow  {\Bbb R}$ is differentiable at any Yamabe metric $\gamma$ 
which is unique in its conformal class $[\gamma]$, for $[\gamma] \neq 
[\gamma_{0}]$, cf. Proposition 2.2. Moreover, if $s$ is differentiable 
at a metric $\gamma$, and $\gamma$ is a critical point of $s$ on 
${\cal Y}^{+}$, then $\gamma$ must be an Einstein metric, (and hence 
unique by Obata's theorem [5]).

 The scalar curvature $s: {\cal Y}^{+} \rightarrow  {\Bbb R}$ is 
continuous, cf. [2]. It is unknown if $s$ is differentiable everywhere; 
at present there are no solid reasons to substantiate this. However, we 
show that away from $[\gamma_{0}]$, positive directional derivatives of 
$s$ always exist, i.e. the derivative 
$$s' (h) = \frac{d}{dt}^{+}s_{[\gamma +th]}|_{t=0}$$
exists, for any symmetric bilinear form $h$; here $s_{[g]}$ denotes the 
scalar curvature of (any) unit volume Yamabe metric in the conformal 
class $[g] \neq [\gamma_{0}]$ and $t > 0$.

 It has long been an open problem whether a Yamabe metric $\gamma$ 
which is a local maximum of $s$ on ${\cal Y}^{+}$ is necessarily an 
Einstein metric, cf. [2]. Some recent progress on this question has 
been made in dimension 3 in [3]. The following result gives a partial 
answer in any dimension.

\begin{theorem} \label{t 1.2.}
  Suppose $\gamma\in{\cal Y}^{+}$ is a local maximum for s. If there 
are most two Yamabe metrics in the conformal class $[\gamma]$, then 
$\gamma $ is Einstein.
\end{theorem}

  I would like to thank the referee for several very useful comments on the paper. 

\section{Proofs of the Results.}

\setcounter{equation}{0}

 Let ${\cal Y}$ denote the space of all unit volume Yamabe metrics on a 
given $n$-manifold $M$. Henceforth, all metrics will be assumed to have 
unit volume, for convenience. The discussion and results to follow are 
all essentially trivial for Yamabe metrics of non-positive scalar 
curvature. Let ${\cal C}$ denote the space of unit volume constant 
scalar curvature metrics on $M$, and ${\Bbb M}$ the space of all unit 
volume metrics on $M$, so that ${\cal Y}  \subset  {\cal C}  \subset  
{\Bbb M}$. Unless stated otherwise, all metrics will be assumed to be 
$(C^{\infty})$ smooth.

 It is well-known, cf. [2], that generically ${\cal C}$ is an infinite 
dimensional manifold. More precisely, at all metrics $g$ for which 
$s_{g}/(n-1)$ is not an eigenvalue of the Laplacian $-\Delta$, (with 
non-negative spectrum), a neighborhood of $g$ in ${\cal C}$ has the 
structure of a smooth manifold. Observe that ${\cal C}$ can be 
described as the 0-level set of the mapping
$$\Phi : {\Bbb M}  \rightarrow  C^{\infty}(M, {\Bbb R} ), \ \Phi (g) = 
\Delta_{g}s_{g}. $$
Since $\Phi$ is a real-analytic map, ${\cal C}$ thus has the structure 
of an infinite dimensional real-analytic variety. It would be 
interesting to know if this structure could lead to further information 
on the structure of ${\cal C}$.

  Consider the set $\Lambda = \Lambda [\gamma]$ of all smooth Yamabe 
metrics in a given conformal class $[\gamma]$, i.e. $\Lambda = {\cal Y} 
\cap [\gamma]$. If $g$ is a fixed representative metric in $\gamma$, 
then all Yamabe metrics $\gamma \in \Lambda [\gamma]$ are of the form
\begin{equation} \label{e2.1}
\gamma = \phi^{4/(n-2)}g,
\end{equation}
for some positive smooth function $\phi: M \rightarrow {\Bbb R}$. By the 
solution to the Yamabe problem [1], [6], the set $\Lambda [\gamma]$ is 
compact, i.e. is a compact subset of $C^{\infty}(M, {\Bbb R})$, for all 
$[\gamma] \neq [\gamma_{0}]$, cf. also [7, Prop. 2.1]. In the following, 
we will always work with conformal classes $[\gamma] \neq [\gamma_{0}]$. 
In addition, the sets $\Lambda [\gamma]$ as $\gamma$ varies are also 
compact in the following sense; if $[\gamma_{i}] \rightarrow  
[\gamma] \neq [\gamma_{0}]$ smoothly, then any sequence of Yamabe metrics 
$(\gamma^{j})_{i} \in  [\gamma_{i}]$ has a subsequence converging smoothly 
to a Yamabe metric $\gamma^{j} \in  [\gamma]$.

\medskip

  We begin by discussing the differentiability properties of the scalar 
curvature function $s$ on ${\cal Y}$. Let $\gamma$ be a unit volume 
Yamabe metric, with $\gamma \notin [\gamma_{0}]$, and let $g_{t}$ be a 
smooth curve of unit volume metrics of the form 
\begin{equation} \label{e2.2}
g_{t} = \gamma + th + O(t^{2}), 
\end{equation}
for $h\in T_{\gamma}{\Bbb M}$. We assume that $g_{t}$ is orthogonal to the 
conformal classes in the sense that
\begin{equation} \label{e2.3}
dV_{g_{t}} = dV_{\gamma}.
\end{equation}
In particular, (2.3) implies that $tr_{\gamma}h = 0$.

  Choose a sequence $t_{i} \rightarrow 0$, and let $\gamma_{i}$ be a unit 
volume Yamabe metric in 
the conformal class of $[g_{t_{i}}]$. Thus, $\gamma_{i}$ is of the form
\begin{equation} \label{e2.4}
\gamma_{i} = \phi_{i}^{4/(n-2)}g_{i}, 
\end{equation}
where $g_{i} = g_{t_{i}}$. By passing to a subsequence if necessary, 
may assume that $\gamma_{i} \rightarrow  \gamma'$ smoothly, where $\gamma'$ 
is some Yamabe metric in $[\gamma]$. Since $g_{i} \rightarrow g$ smoothly, 
one has $\phi_{i} \rightarrow \phi$ smoothly, where $\gamma'  = 
\phi^{4/(n-2)}\gamma$. In general, $\phi  \neq 1$, since the conformal class 
$[\gamma]$ may have more than one Yamabe metric. However, by changing the curve 
$g_{t}$ to the curve $\widetilde g_{t} = \phi^{4/(n-2)}g_{t}$, we may and do assume, 
without loss of generality, that $\phi  = 1$. This simplifies some of the computations 
below. 

 The functions $\phi_{i}$ satisfy the Yamabe equation
\begin{equation} \label{e2.5}
-c_{n}\Delta_{g_{i}}\phi_{i} + s_{g_{i}}\phi_{i} = 
\phi_{i}^{q}s_{\gamma_{i}}, 
\end{equation}
where $c_{n} = 4(n-1)/(n-2)$ and $q = (n+2)/(n-2)$.

 Next, linearize (2.5) about the limit $(\gamma, 1)$ of $(g_{i}, 
\phi_{i})$. The linearization of the scalar curvature is given by 
$$L(h) = \frac{d}{dt}s_{g+th}|_{t=0} = \Delta tr h + \delta\delta h - 
\langle Ric, h \rangle, $$
cf. [2] for instance. The adjoint operator $L^{*}$ to $L$ is given by
$$L^{*}f = D^{2}f - \Delta f\cdot  g - fRic. $$
Since $g_{i} \rightarrow \gamma$ smoothly,
\begin{equation} \label{e2.6}
s_{g_{i}} = s_{\gamma} + t_{i}L(h) + O(t_{i}^{2}).
\end{equation}
Also write
$$\phi_{i}^{q}s_{\gamma_{i}} = (\phi_{i}^{q}-\phi_{i})s_{\gamma} + 
\phi_{i}^{q}(s_{\gamma_{i}}-s_{\gamma}) + \phi_{i}s_{\gamma}.$$
Substituting these expressions in (2.5) gives
$$-c_{n}\Delta_{g_{i}}(\phi_{i}-1) - (\phi_{i}^{q}-\phi_{i})s_{\gamma} 
= \phi_{i}^{q}(\frac{s_{\gamma_{i}}-s_{\gamma}}{t_{i}})t_{i} - 
\phi_{i}L(h)t_{i} + O(t_{i}^{2}), $$
so that
\begin{equation} \label{e2.7}
-c_{n}\Delta_{g_{i}}(\frac{\phi_{i}-1}{t_{i}}) - 
\psi_{i}s_{\gamma}(\frac{\phi_{i}-1}{t_{i}}) = 
\phi_{i}^{q}(\frac{s_{\gamma_{i}}-s_{\gamma}}{t_{i}}) -\phi_{i}L(h)+O(t_{i}),
\end{equation}
where
\begin{equation} \label{e2.8}
\psi_{i} = (\phi_{i}^{4}+...+\phi_{i})/((\phi_{i}^{q-1})^{n-3}+ 
(\phi_{i}^{q-1})^{n-4}+ ...+\phi_{i}^{q-1}+1). 
\end{equation}
Here we use the identity that if $y = x^{a/b}$, then $(y-1) = 
(x-1)(x^{a-1}+...+1)/(y^{b-1}+...+1)$, with $x = \phi$, $y = 
\phi^{q-1}$, $a = 4$, $b = n-2$.

 Integrating (2.7) over $M$ with respect to $g_{i}$ then gives the identity
\begin{equation} \label{e2.9}
-s_{\gamma} \int (\frac{\phi_{i}-1}{t_{i}})\psi_{i}dV_{g_{i}} = 
(\frac{s_{\gamma_{i}}-s_{\gamma}}{t_{i}})\int \phi_{i}^{q}dV_{g_{i}} - 
\int \phi_{i}L(h) dV_{g_{i}} + O(t_{i}). 
\end{equation}

 To understand the linearization, i.e. the behavior as $t_{i} \rightarrow 0$, 
the main point is to prove that the left side of (2.9) vanishes in 
the limit $t_{i} \rightarrow 0$.

\begin{proposition} \label{p 2.1.}
  Let $g_{t}$ be as in (2.2), with $t = t_{i} \rightarrow 0$, $t_{i} > 0$, 
and $s_{\gamma} > 0$. Let $\gamma_{i}$ be Yamabe metrics in $[g_{i}]$, with 
$\phi_{i}$ as in (2.4) satisfying $\phi_{i} \rightarrow 1$ smoothly as 
$t_{i} \rightarrow 0$. Then 
\begin{equation} \label{e2.10}
\int (\frac{\phi_{i}-1}{t_{i}})\psi_{i}dV_{g_{i}} \rightarrow  0, \ 
{\rm as} \  i \rightarrow  \infty . 
\end{equation}
\end{proposition}
\noindent
{\bf Proof:} Note first that since $\phi_{i}^{q} \rightarrow 1$ and 
$\psi_{i} \rightarrow \psi = \frac{4}{n-2}$ smoothly, the only terms 
in (2.9) which may become unbounded are those with $t_{i}$ in the 
denominator. Also, the measures $dV_{g_{i}}$ converge smoothly to the 
measure $dV_{g}$. 

   Suppose first that $(\phi_{i}-1)/t_{i}$ is bounded in $L^{2}$, so that 
it has a weakly convergent subsequence: $(\phi_{i}-1)/t_{i} \rightarrow 
\phi'$ weakly in $L^{2}$. The left side of (2.9) is then bounded, and hence 
the term $(s_{\gamma_{i}}-s_{\gamma})/t_{i}$ on the right is also bounded. 
It then follows from elliptic regularity associated to the equation (2.7) 
that $(\phi_{i}-1)/t_{i} \rightarrow \phi'$ smoothly. Thus, it suffices 
to show in this case that
\begin{equation} \label{e2.11}
\int\phi' dV_{\gamma} = 0, 
\end{equation}
since $\psi  = \lim \psi_{i} = 4/(n-2)$. However, since the metrics 
$s_{g_{i}}$ and $s_{\gamma_{i}}$ have unit volume, one has
\begin{equation} \label{e2.12}
\int\frac{\phi_{i}^{2n/(n-2)} - 1}{t_{i}}dV_{g_{i}} = 0. 
\end{equation}
Taking the limit of (2.12) as $t_{i} \rightarrow 0$ and using the fact 
that $\phi'$ exists gives (2.11).

  If $(\phi_{i}-1)/t_{i}$ is not bounded in $L^{2}$, the proof is more 
complicated, but based on similar ideas together with the fact that 
$\gamma_{i}$ are Yamabe metrics. To begin, we assume that 
$(s_{\gamma_{i}}-s_{\gamma})/t_{i}$ is bounded, and hence converges 
to a limit
\begin{equation} \label{e2.13}
s' = \lim_{t_{i} \rightarrow 0} \frac{s_{\gamma_{i}}-s_{\gamma}}{t_{i}},
\end{equation}
(again in a subsequence). The situation where (2.13) does not hold is 
dealt with later, based on a simple renormalization argument. 
The assumption (2.13) implies that the limit of the left side of (2.9) as 
$t_{i} \rightarrow 0$ also exists.

  First, we claim that
\begin{equation} \label{e2.14}
\lim_{t_{i} \rightarrow 0} -\int(\frac{\phi_{i}-1}{t_{i}})\psi_{i}dV_{g_{i}} 
\leq 0.
\end{equation}
To see this, since $\phi_{i} \rightarrow 1$ smoothly, one has
$$\lim_{t_{i} \rightarrow 0} \frac{s_{\gamma_{i}}-s_{\gamma}}{t_{i}}\int \phi_{i}^{q}dV_{g_{i}} = s'.$$
On the other hand, by the Yamabe property of $\gamma_{i}$ and the fact that 
$g_{i}$ is of unit volume, 
\begin{equation} \label{e2.15}
s_{\gamma_{i}} \leq \int s_{g_{i}}dV_{g_{i}}.
\end{equation}
Hence, by (2.13),
$$s' \leq \lim_{i \rightarrow \infty} \int \frac{s_{g_{i}}-s_{\gamma}}{t_{i}}
dV_{g_{i}}.$$
However, (2.6) gives
$$\lim_{t_{i} \rightarrow 0}\int \frac{s_{g_{i}}-s_{\gamma}}{t_{i}}dV_{g_{i}} 
= \int L(h)dV_{\gamma}.$$
Via (2.9), this gives the claim (2.14).

  We now claim the opposite inequality to (2.14) holds. To see this, consider 
the sequence of metrics $\widetilde g_{i} = \phi_{i}^{4/(n-2)}\gamma$ in the 
conformal class $[\gamma]$. They all have unit volume, by (2.2) and (2.12). 
Since $\gamma$ is a Yamabe metric, i.e. it minimizes ${\cal S}$ in its 
conformal class, it follows that
\begin{equation} \label{e2.16}
\int (c_{n}|d\phi_{i}|^{2} + s_{\gamma}\phi_{i}^{2})dV_{\gamma} = 
\int [c_{n}|d(\phi_{i} - 1)|^{2} + s_{\gamma}(1 + (\phi_{i} - 1))^{2}]
dV_{\gamma} \geq s_{\gamma}.
\end{equation}
Expanding out the term on the right then gives, since $t_{i} > 0$ and 
$s_{\gamma} > 0$,
$$\int [c_{n}\frac{|d(\phi_{i} - 1)|^{2}}{t_{i}} + 
s_{\gamma}(2\frac{(\phi_{i}-1)}{t_{i}} +\frac{(\phi_{i}-1)^{2}}{t_{i}})] 
dV_{\gamma} \geq 0.$$
Also, integration by parts gives
\begin{equation} \label{e2.17}
\int c_{n}\frac{|d(\phi_{i} - 1)|^{2}}{t_{i}}dV_{\gamma} = 
-c_{n}\int (\phi_{i}-1)\Delta_{\gamma}({\frac{\phi_{i}-1}{t_{i}}})dV_{\gamma} 
= s_{\gamma}\int \psi_{i}\frac{(\phi_{i}-1)^{2}}{t_{i}}dV_{\gamma} + o(1).
\end{equation}
Here the second equality uses (2.7), together with the fact that 
$\frac{1}{t_{i}}(\Delta_{g_{i}}- \Delta_{\gamma})(\phi_{i}-1) \rightarrow 0$ 
since $g_{i} \rightarrow \gamma$ smoothly. (The term 
$\frac{1}{t_{i}}(\Delta_{g_{i}}- \Delta_{\gamma})$ converges to the derivative 
$(\Delta)'_{h}$ of the Laplacian in the direction $h$; this applied to 
$(\phi_{i} - 1)$ tends to 0, since $(\phi_{i} - 1) \rightarrow 0$). Combining 
these estimates, (and using $s_{\gamma} > 0$), it follows that 
\begin{equation} \label{e2.18}
\lim_{t_{i} \rightarrow 0}\int 2\frac{(\phi_{i}-1)}{t_{i}}+\frac{(\phi_{i}-1)^{2}}{t_{i}} + 
\psi_{i}\frac{(\phi_{i}-1)^{2}}{t_{i}}dV_{\gamma} \geq 0.
\end{equation}
Observe that from the derivation of (2.7), we have 
$\psi_{i}(\phi_{i}-1)^{2} = (\phi_{i}-1)(\phi_{i}^{q}-\phi_{i})$. Now 
compute $2(\phi_{i}-1)+(\phi_{i}-1)^{2}+(\phi_{i}-1)(\phi_{i}^{q}-\phi_{i}) 
= (\phi_{i}-1)(2+(\phi_{i}-1)+ \phi_{i}^{q}-\phi_{i}) = (1+\phi_{i}^{q})
(\phi_{i}-1) = \phi_{i}-\phi_{i}^{q} + \phi_{i}^{q+1} - 1$. The last two 
terms here integrate to 0, by (2.12), since $q+1 = 2n/(n-2)$. It follows 
then from (2.18) that
$$\lim_{t_{i}\rightarrow 0}-\int \frac{\phi_{i}^{q} - \phi_{i}}{t_{i}}dV_{\gamma} 
\geq 0,$$
i.e.
\begin{equation} \label{e2.19}
\lim_{t_{i}\rightarrow 0}-\int \psi_{i}\frac{\phi_{i}-1}{t_{i}} \geq 0.
\end{equation}
Combining (2.14) and (2.19) proves the result in case (2.13) holds.

  Finally, suppose the ratio $(s_{\gamma_{i}} - s_{\gamma})/t_{i}$ in (2.13) 
is unbounded. Then divide each term in (2.9) by $C_{i}$, where $C_{i} \rightarrow 
\infty$ is chosen to make the resulting ratio in (2.13) equal to 1, in absolute 
value; thus $(s_{\gamma_{i}} - s_{\gamma})/C_{i}t_{i}$ remains bounded, and bounded 
away from 0. Performing exactly the arguments as above following (2.13), dividing by 
$C_{i}$ as called for, leads to the conclusion that $\int (\phi_{i}-1)/C_{i}t_{i} 
\rightarrow 0$, and $\int \phi_{i}L(h)/C_{i} \rightarrow 0$, but $(s_{\gamma_{i}} - 
s_{\gamma})/C_{i}t_{i}$ is bounded away from 0. This contradicts (2.9), and so 
completes the proof. 
{\endproof}

 Proposition 2.1 and (2.9) imply that the limit
\begin{equation} \label{e2.20}
s'(h) = \lim_{t_{i}\rightarrow 0^{+}}\frac{s_{\gamma_{i}}-s_{\gamma}}{t_{i}},
\end{equation} 
exists, and is given by
\begin{equation} \label{e2.21}
s' (h) = \int \langle L^{*}1, h \rangle dV_{\gamma} =\int \langle -z, h 
\rangle dV_{\gamma}, 
\end{equation}
where $z$ is the trace-free Ricci curvature, $z = Ric - \frac{s}{n}g$. 

 The discussion above easily proves the following result.

\begin{proposition} \label{p 2.2.}
  Suppose $\gamma $ is the unique Yamabe metric (of unit volume) in 
$[\gamma] \neq [\gamma_{0}]$. Then $s$ is differentiable on ${\cal Y}$ 
at $\gamma$, and the derivative is given by (2.21).
\end{proposition}
\noindent
{\bf Proof:} Let $[g_{t}] = [\gamma +th]$ be any variation of the 
conformal class of $\gamma$. If $\gamma_{t}$ is any Yamabe metric in 
$[g_{t}]$ written in the form (2.4), then as discussed following (2.4), 
$\phi_{t} \rightarrow 1$ as $t \rightarrow 0$. It follows from (2.21) 
that for any $h$, and for any sequence of Yamabe metrics converging to 
$\gamma$ along the curve $[g_{t}]$, one has
$$s' (h) = \int \langle -z, h \rangle dV_{\gamma}. $$
The right hand side is linear in $h$, and hence $s$ is differentiable 
at $\gamma .$
{\endproof}

\begin{corollary} \label{c 2.3.}
  Any critical point $\gamma$ of $s$ on ${\cal Y}$, for which $\gamma$ 
is unique in its conformal class is necessarily an Einstein metric. 
\end{corollary}
\noindent
{\bf Proof:} If $\gamma$ is unique in its conformal class, then 
Proposition 2.2 implies that $s$ is differentiable on ${\cal Y}$ at 
$\gamma$. By (2.21), a critical point $\gamma$ of $s$ then satisfies $z 
= 0$, i.e. $\gamma$ is Einstein.

{\endproof}

 Now suppose there is more than one Yamabe metric in the conformal 
class $[\gamma]$. Such metrics are of the form $\widetilde \gamma = 
\phi^{4/(n-2)}\gamma$, for some positive function $\phi$. A standard 
formula for conformal changes, cf. [2] gives, 
$$\widetilde z = z + (n-2)\phi^{2/(n-2)}D_{0}^{2}\phi^{-2/(n-2)}, $$
where $D_{0}^{2}$ is the trace-free Hessian with respect to $\gamma$. 
Since $h$ transforms as $\widetilde h = \phi^{4/(n-2)}h$ on passing 
from $\gamma$ to $\widetilde \gamma$, and $dV_{\widetilde \gamma} = 
\phi^{2n/(n-2)}dV_{\gamma}$, one has
$$\int \langle \widetilde z, \widetilde h \rangle dV_{\widetilde \gamma} =  
\int \phi^{(2n-4)/(n-2)} \langle z_{\gamma} + 
(n-2)\phi^{2/(n-2)}D_{0}^{2}\phi^{-2/(n-2)} , h \rangle dV_{\gamma}.$$
Set $u = \phi^{2/(n-2)}$, so that
\begin{equation} \label{e2.22}
\int \langle \widetilde z, \widetilde h \rangle dV_{\widetilde \gamma} =  
\int u^{n-2} \langle z_{\gamma} + (n-2)u(D_{0}^{2}u^{-1}) , h \rangle 
dV_{\gamma}.
\end{equation}
The forms
\begin{equation} \label{e2.23}
Z_{\phi} = u^{n-2}(z_{\gamma} + (n-2)uD_{0}^{2}u^{-1}), 
\end{equation}
are viewed as tangent vectors $Z_{\phi}\in T_{\gamma}{\Bbb M}$, 
corresponding to the set of Yamabe metrics $\Lambda [\gamma]$ in $[\gamma]$, 
cf. the discussion preceding (2.1). Let $g_{t}$, $t \geq 0$, be a smooth curve 
of unit volume metrics as in (2.2) and let $\gamma_{i}$ be an associated sequence 
of Yamabe metrics as in (2.4). It follows from (2.21) and (2.22)-(2.23) that if 
$\gamma_{i}$ converges to a Yamabe metric 
$\widetilde \gamma = \phi^{4/(n-2)}\gamma$ in $[\gamma] = [g(0)]$, then
\begin{equation} \label{e2.24}
s' (h) = \int \langle -Z_{\phi}, h \rangle dV_{\gamma}, 
\end{equation}
where $s' = \lim_{t_{i} \rightarrow 0^{+}}(s_{\gamma_{i}}-s_{\gamma})/t_{i}$. 

  One thus has a map
\begin{equation} \label{e2.25}
{\cal Z}: \Lambda [\gamma] \rightarrow  \{Z_{\phi}\}, \ \  
{\cal Z}(\phi^{4/(n-2)}\gamma) = Z_{\phi}. 
\end{equation}
It will be shown below that ${\cal Z}$ is 1-1, cf. Lemma 2.5. While of 
course $Z_{\phi}$ {\it may} be defined for any Yamabe metric in $\Lambda [\gamma]$, 
we only define $Z_{\phi}$ for a Yamabe metric $\widetilde \gamma$ in $[\gamma]$ 
arising as a limit of a sequence $\gamma_{i}$ as in (2.24), and call such $Z_{\phi}$ 
admissible. It is not known if all $Z_{\phi}$ are admissible, i.e. if ${\cal Z}$ is 
then defined on all of $\Lambda [\gamma]$. There may Yamabe metrics $\bar \gamma$ in 
$[\gamma]$ which are not moved in any perturbation of the conformal class 
$[\gamma]$; such $\bar \gamma$ are thus not near a Yamabe metric in any conformal 
class close to $[\gamma]$, and so are ``isolated'' Yamabe metrics. 

  It is a consequence of (2.24) that $s'$ is independent of the 
sequence $t_{i}$, among sequences for which $\gamma_{i} \rightarrow \widetilde 
\gamma$ as above. However, apriori, even within the fixed curve $g_{t}$ specified in 
(2.24), different sequences $t_{i} \rightarrow 0$ may give sequences of Yamabe metrics 
converging to different Yamabe limits in $[\gamma]$. Thus, apriori, one 
may have a collection of distinct $\phi$'s associated to a given $h$. 
Of course in general the collection of $\phi$'s and $Z_{\phi}$'s 
changes with $h$. 

 This analysis leads to the following result.

\begin{proposition} \label{p 2.4.}
  The scalar curvature function s: ${\cal Y}  \rightarrow  {\Bbb R}$ 
has (positive) directional derivatives in any direction 
$h\in T_{\gamma}{\Bbb M}$, $[\gamma] \neq [\gamma_{0}]$, with derivative 
given by
\begin{equation} \label{e2.26}
s' (h) = \frac{d}{dt}^{+}s_{[\gamma +th]}|_{t=0} = \min_{\phi}\int \langle 
-Z_{\phi}, h \rangle dV_{\gamma}, 
\end{equation}
where the minimum is taken over all admissible $Z_{\phi}$ in $[\gamma]$.
\end{proposition}
\noindent
{\bf Proof:}
This is an immediate consequence of the discussion above and the fact 
that Yamabe metrics minimize the scalar curvature functional ${\cal S}$ 
in the conformal class; in particular, $s_{\gamma}$ is independent of 
the representative $\gamma\in [\gamma]$.
{\endproof}

 As above, let $\{\gamma^{\lambda}\}_{\lambda\in\Lambda}$ denote the set of 
unit volume Yamabe metrics in the conformal class $[\gamma]$, 
$[\gamma] \neq [\gamma_{0}]$. By the compactness mentioned at the beginning 
of this section, if $[\gamma_{i}] \rightarrow [\gamma]$ smoothly, then any 
sequence of Yamabe metrics $(\gamma^{j})_{i} \in [\gamma_{i}]$ has a subsequence 
converging smoothly to some Yamabe metric $\gamma^{j} \in [\gamma]$. However, 
the cardinality of the sets $\{\gamma^{\lambda}\}$ may well change in passing 
to limits. The cardinality may drop, for instance when distinct Yamabe metrics 
in a sequence of conformal classes merge to a common limit. The cardinality 
of $\{\gamma^{\lambda}\}$ may also increase in the limit, due to the 
``birth'' of a new Yamabe metric, not arising as a limit of a given 
sequence. 

\medskip

 We are now in position to prove Theorem 1.1.

{\bf Proof of Theorem 1.1.}

 It is well-known, cf. [8] for instance, that the set of conformal 
classes containing only finitely many unit volume Yamabe (or more 
generally constant scalar curvature) metrics is generic, i.e. forms an 
open and dense set in the space of all Yamabe metrics, with respect to the 
$C^{2}$ (or stronger) topology. For simplicity, we assume $[\gamma]$ is 
generic in this sense. Then nearby conformal classes $[g_{t}]$, in any 
direction $g_{t} = \gamma +th$, are also generic, and so have only finitely 
many Yamabe metrics.

 The proof is by contradiction, and so suppose there exists an open set 
${\cal V}$ of generic Yamabe metrics such that, for each $\gamma\in{\cal V}$ 
there are at least two distinct Yamabe metrics in $[\gamma]$. Using the 
compactness mentioned above, by passing to a smaller ${\cal V}$ if necessary, 
we may assume that there is a fixed upper bound on the number of distinct Yamabe 
metrics in $[\gamma]$, $\gamma\in{\cal V}$.

 Now fix some $\gamma_{o}\in{\cal V}$, and consider the class of 
variations $h\in T_{\gamma_{o}}{\Bbb M}$ with $||h||_{C^{2}} = 1$. As 
discussed above, each such $h$ determines a subset of Yamabe metrics 
$\{\gamma^{\lambda}(h)\}$ of the full set of Yamabe metrics 
$\{\gamma^{\lambda}\}$ in $[\gamma_{o}]$, (namely those which persist 
under perturbation in the direction $h$). As $h$ changes, this subset 
can of course change. However, since there is a fixed bound on the 
cardinality of $\{\gamma^{\lambda}\}$, there is an open set $U$ of 
$h\in T_{\gamma}{\Bbb M} \cap \{||h||_{C^{2}} = 1\}$, with the following 
property: there is a fixed subset $\{\gamma^{\lambda'}\} \subset 
\{\gamma^{\lambda}\}$, of cardinality at least two, such that for each 
$h\in U$, there is a sequence $t_{i} \rightarrow 0^{+}$ such that 
$[\gamma_{o}+t_{i}h]$ has at least two Yamabe metrics converging in 
subsequences to at least two distinct elements of $\{\gamma^{\lambda'}\}$. 
We may assume that $\gamma_{o}\in\{\gamma^{\lambda'}\}$. Henceforth, 
set $\gamma_{o} = \gamma$ to simplify the notation. 

 Since the scalar curvature of a Yamabe metric depends only on the conformal 
class of the metric, it then follows that, for all $h\in U$, one has
\begin{equation} \label{e2.27}
s'(h) = \int \langle -Z_{\phi}, h \rangle dV_{\gamma} = \int \langle -z, h 
\rangle dV_{\gamma}, 
\end{equation}
for all $\phi$ representing the Yamabe metrics in 
$\{\gamma^{\lambda'}\}$. In particular, (2.27) holds for at least one 
$\phi  \neq  1$. The condition (2.27) is linear and since it holds for 
all $h\in U$, it must hold for all $h$ in the linear span of $U$, and 
hence it holds for all (trace-free) $h\in T_{\gamma_{o}}{\Bbb M}$. This 
implies that one has pointwise
\begin{equation} \label{e2.28}
Z_{\phi} = z, 
\end{equation}
i.e. 
\begin{equation} \label{e2.29}
u^{n-2}(z + (n-2)uD_{0}^{2}u^{-1}) = z, 
\end{equation}
for some $\phi \neq 1$. 

 Thus to prove Theorem 1.1, it suffices to prove that (2.29) has only 
trivial solutions. This is done by a computation in the following Lemma.
\begin{lemma} \label{l 2.5.}
  For $(M, [\gamma]) \neq (S^{n}, [\gamma_{0}])$, the only solution to (2.29) 
is 
\begin{equation} \label{e2.30}
u = \phi  = 1.
\end{equation}
\end{lemma}
\noindent
{\bf Proof:}
We first prove this in case $n > 3$ and afterwards prove the case $n = 3$. 
Write (2.29) as
\begin{equation} \label{e2.31}
u^{-(n-2)}(1-u^{n-2})z = (n-2)u(D_{0}^{2}u^{-1}), 
\end{equation}
where as above, $u = \phi^{2/(n-2)}$. For $\widetilde g = u^{2}g$, as noted 
following Corollary 2.3, 
$$\widetilde z = z + (n-2)uD_{0}^{2}u^{-1}.$$
Combining this with (2.31) gives
\begin{equation} \label{e2.32}
z = \phi^{2}\widetilde z.
\end{equation}
The metrics $\widetilde g$ and $g$ are Yamabe, and so by the Bianchi identity, 
$$\delta_{g}(\phi^{2}\widetilde z) = \phi^{2}\delta_{g}\widetilde z - 
\widetilde z(\nabla_{g}\phi^{2}) = 0.$$
A simple computation shows that 
$$\delta_{g}\widetilde z = u^{2}\delta_{\widetilde g}\widetilde z - 
u^{2}\widetilde z(\nabla_{\widetilde g}\log u^{-1}).$$
Since $\nabla_{g}\phi^{2} = u^{2}\nabla_{\widetilde g}\phi^{2}$, 
these equations and (2.32) give
$$\widetilde z (\nabla_{\widetilde g}\log u) = 
\widetilde z (\nabla_{\widetilde g}\log \phi^{2}).$$
Since $\phi^{2} = u^{n-2}$, it follows that $\widetilde z(\nabla_{\widetilde g}u) = 
0$, because $n > 3$. Interchanging the roles of $g$ and $\widetilde g$ then gives 
\begin{equation} \label{e2.33}
z(\nabla_{g}u) = 0.
\end{equation}

  To complete the proof when $n > 3$, return to (2.31). Let $f(u) = 
u^{-(n-1)}(1-u^{n-2})$. Assuming $f(u)$ is not identically $0$, divide 
(2.31) by $f(u)$ and pair it with $D^{2}u^{-1}$ to obtain
\begin{equation} \label{e2.34}
\langle D^{2}u^{-1}, z \rangle = \frac{(n-2)}{f(u)} \langle D^{2}u^{-1}, 
D^{2}u^{-1} - \frac{\Delta u^{-1}}{n}g \rangle.
\end{equation}
The left side of (2.34) is smooth, and hence so is the right side. Let 
$U^{+} = \{u \geq 1 \}$. Applying the divergence theorem over the domain 
$U^{+}$ gives
$$\int_{U^{+}} \langle D^{2}u^{-1}, z \rangle = \int_{U^{+}} 
\langle du^{-1}, \delta z \rangle + \int_{\partial U^{+}}z(du^{-1}, \nu) 
= 0,$$
where we have used the Bianchi identity and (2.33). Hence
$$\int_{U^{+}} \frac{n-2}{f(u)} \{|D^{2}u^{-1}|^{2} - 
\frac{1}{n}(\Delta u^{-1})^{2} \} = 0.$$
Since $f(u) \leq 0$ in $U^{+}$, the Cauchy-Schwarz inequality shows that the integrand 
here is pointwise non-positive, and so one must have 
\begin{equation} \label{e2.35}
D^{2}u^{-1} = \frac{\Delta u^{-1}}{n}g, 
\end{equation}
on $(U^{+}, \gamma)$. The same arguments apply to the complementary region $U^{-}$, 
and so (2.35) holds on $(M, \gamma)$. By (2.31), this of course implies $z = 0$, at 
least on the domain where $f(u) \neq 0$. Taking the divergence of (2.35) on this 
domain and using the identity $div D^{2} \psi = d \Delta \psi + Ric(\nabla \psi)$ 
gives
$$d(\frac{n-1}{n} \Delta u^{-1} + \frac{su^{-1}}{n})= 0.$$
Since $u$ is smooth on $M$, it follows that $\Delta u^{-1} + \frac{s}{n-1}u^{-1} = c$ 
on $M$, with $c = \frac{s}{n-1}\int u^{-1}$. Let $v = u^{-1} - \frac{n-1}{s}c$, so 
that $\int v = 0$, and $\Delta v = -\frac{s}{n-1}v$. It follows then from (2.35) that 
\begin{equation} \label{e2.36}
D^{2}v = -\frac{s}{n(n-1)}vg.
\end{equation}
By Obata's theorem [5], it is well-known that the only non-zero solution to (2.36) 
is with $v$ a $1^{\rm st}$ eigenfunction of the Laplacian on the round metric 
$\gamma_{0}$ on $S^{n}$, (up to a constant rescaling). Since by assumption 
$(M, [\gamma]) \neq (S^{n}, [\gamma_{0}])$, it follows 
then that $v = 0$ and hence $u = 1$. This completes the proof in case $n > 3$. 

  For $n = 3$, a different argument must be used, since the proof of (2.33) does not 
hold in this case. Instead, return to (2.32), and recall that $\omega = 
d(r - \frac{s}{4}g)$ is a conformal invariant of weight 0 in dimension 3, cf. [2]; 
thus if $\widetilde g = u^{2}g$, then $\widetilde \omega = \omega$. The metrics $g$ 
and $\widetilde g$ are Yamabe, and so 
$$\widetilde{dz} = dz.$$
On the other hand, by (2.32), one has 
$$dz = d\phi^{2} \wedge \widetilde z + \phi^{2}d\widetilde z,$$
and one easily computes from the definition that $d\widetilde z = \widetilde{dz}$. 
Hence $(1 - \phi^{2})dz = d\phi^{2} \wedge z$. This implies that 
$d\phi^{2} \wedge z = 0$ on the level set $\{\phi = u = 1\}$. Taking the trace 
of this $(2,1)$ form gives $z(du) = 0$ on $\{u = 1\}$. The remainder of the proof 
then follows exactly as following (2.34). 
{\endproof}

  The proof of Theorem 1.1 now follows easily. Lemma 2.5 implies there 
are no non-trivial solutions of (2.29). Hence the open set ${\cal V}$ 
must be empty, which proves the result.

{\endproof}

 The proof of Theorem 1.2 is now also very simple.

{\bf Proof of Theorem 1.2.}

 Suppose that the scalar curvature function $s$ has a local maximum on 
${\cal Y}$ at $\gamma$, i.e. $s_{\gamma'} \leq  s_{\gamma}$, for any 
unit volume Yamabe metric $\gamma'$ near $\gamma$. Then by Proposition 
2.4, the positive directional derivative $s'(h)$ exists for any $h$ and
$$s' (h) \leq  0. $$
Consider then the collection of tangent vectors $Z_{\phi}$ in 
$T_{\gamma}{\Bbb M}$ with $Z_{\phi} = z$ for $\phi = 1$. Any $h$ determines 
at least one $Z_{\phi}$, (and maybe several), and (2.26) gives
$$\int \langle -Z_{\phi}, h \rangle dV_{\gamma} \leq  0. $$
Hence, each $h$ must lie in the positive half-space of its 
corresponding $Z_{\phi}$. However, if there are only two Yamabe metrics 
in $[\gamma]$, there are most two such $Z_{\phi}$ and the two 
corresponding half-spaces do not cover all of $T_{\gamma}{\Bbb M}$ 
unless
$$z = -Z_{\phi}.$$
From (2.23), this implies
\begin{equation} \label{e2.37}
u^{-1}(1+u^{n-2})z = -(n-2)D_{0}^{2}u^{-1}. 
\end{equation}
We now argue as in the proof of Lemma 2.5. Thus, pair (2.37) with $z$; integrating 
over $M$ and using the Bianchi identity implies that $z = 0$, i.e $\gamma$ is Einstein.
{\endproof}

 It would be of interest to know if this method of proof can be 
generalized to prove that any local maximum of $s$ on ${\cal Y}$ is 
necessarily Einstein.

\bibliographystyle{plain}

\bigskip

\begin{center}
May, 2003
\end{center}

\medskip
\noindent
\address{Department of Mathematics\\
S.U.N.Y. at Stony Brook\\
Stony Brook, N.Y. 11794-3651\\}
\noindent
\email{anderson@@math.sunysb.edu}

\end{document}